\documentclass[10pt]{article}
\usepackage{amssymb}
\usepackage{amsmath}
\usepackage{amsbsy}
\newtheorem{lemma}{Lemma}[section]
\newtheorem{proposition}{Proposition}[section]
\newtheorem{theorem}{Theorem}[section]

\newtheorem{remark}{Remark}[section]

\def\proclaim#1{\par \bigskip\noindent {\bf #1}\bgroup\it\ }
\def\endproclaim{\egroup\par\bigskip}

\newbox\TempBox \newbox\TempBoxA
\newcommand{\non}{\nonumber \\}
\def\pr{\textsf{P}} % the symbol P for probability used the sans serif letter
\def\ep{\textsf{E}} % the symbol E for expectation used the sans serif letter
\def\Var{\textsf{Var}} % the symbol Var for covariance used the sans serif letter
\def\Cal#1{{\mathcal #1}}

\def\text#1{\mbox{\rm #1}}
\def\overset#1#2{\stackrel{#1}{#2} }

\def\underwiggle 1{
\ifmmode\setbox\TempBox=\hbox{$ 1$}\else\setbox\TempBox=\hbox{
1}\fi \setbox\TempBoxA=\hbox to \wd\TempBox{\hss\char'176\hss}
\rlap{\copy\TempBox}\smash{\lower9pt\hbox{\copy\TempBoxA}} }

\begin{document}

%\renewcommand{\theequation}
%{\arabic{section}.\arabic{equation}}
%\baselineskip=22pt

%$ $ \vskip 0.8in

%\Large

\title{\huge \bf Precise rates in the law of the iterated logarithm
\footnote{Research supported by  National Natural Science Foundation
of China }}
\author{
{\sc   Li-Xin Zhang \footnote{Department of Mathematics, Zhejiang
University, Hangzhou 310027, China,\newline E-mail:
lxzhang@stazlx.zju.edu.cn}}
\\
{\em Department of Mathematics, Zhejiang University, Hangzhou
310027, China }}

%\date{\today}
\date{ }

\maketitle

{\rm
%\sf

{\sc Abstract.} \quad Let $X$, $X_1$, $X_2$, $\ldots$ be i.i.d.
random variables, and set $S_n=X_1+\ldots + X_n$, $M_n=\max_{k\le
n}|S_k|$, $n\ge 1$. Let $a_n=o(\sqrt{n/\log\log n})$. By using the
strong approximation, we prove that,  if $\ep X^2I\{|X|\ge
t\}=o((\log\log t)^{-1})$ as $t\to \infty$, then for $a>-1$ and
$b>-1$,
\begin{eqnarray*}
& & \lim_{\epsilon\searrow \sqrt{1+a}}
(\epsilon^2-a-1)^{b+1/2} \sum_{n=1}^{\infty}
\frac{(\log n)^a (\log\log n)^b}{n}
\pr \Big\{M_n \ge \sigma \phi(n) \epsilon
+a_n\Big\}
\\
& &= 2\sqrt{\frac 1{\pi (a+1)} } \Gamma(b+1/2)
\end{eqnarray*}
holds if and only if
$$ \ep X= 0,\;\; \ep X^2=\sigma^2<\infty \;
\text{ and  }\;
\ep\big[ X^2(\log |X|)^a(\log\log |X|)^{b-1}\big]<\infty.$$
We also show that the condition $\ep X^2I\{|X|\ge t\}=o((\log\log t)^{-1})$ is sharp. The results of Gut and Sp\u ataru (2000) are special cases of ours.

\bigskip

{\bf Keywords:} Tail probabilities of sums of i.i.d. random variables,
\quad the law of the iterated logarithm,  \quad strong approximation.

\bigskip
{\bf AMS 1991 subject classification:} Primary 60F15, Secondary
60G50.

 \vskip 0.2in

}
\newpage

%\newpage

\section{Introduction and main results.}
\setcounter{equation}{0}

Let $\{X, X_n; n\ge 1\}$ be a sequence of i.i.d random variables
with common distribution function $F$, mean $0$ and positive,
finite variance $\sigma^2$, and set $S_n=\sum_{k=1}^n X_k$,
$M_n=\max_{k\le n}|S_k|$, $n\ge 1$. Also let $\log x=\ln(x\vee
e)$, $\log\log x=\log(\log x)$ and $\phi(x)=\sqrt{2x\log\log x}$.
  Then by the well-known law of the iterated logarithm
  (LIL) we have
  \begin{equation} \label{eq1.1}
  \limsup_{n\to \infty}\frac{M_n}{\phi(n)}
=\limsup_{n\to \infty}\frac{|S_n|}{\phi(n)}
=\sigma \quad a.s..
  \end{equation}
Gut and Sp\u ataru (2000) proved the
following two results on its precise asymptotics.

\proclaim{Theorem A} Suppose that $\ep X=0$, $\ep X^2=\sigma^2$ and $\ep[X^2(\log\log |X|)^{1+\delta}]<\infty$
for some $\delta>0$, and let
$a_n=O(\sqrt{n}/ (\log\log n)^{\gamma})$ for some $\gamma>1/2$. Then
$$\lim_{\epsilon\searrow 1}\sqrt{\epsilon^2-1}\sum_{n=1}^{\infty}\frac 1n
\pr(|S_n|\ge \epsilon \sigma \phi(n)+a_n)=1 . $$
\endproclaim

\proclaim{Theorem B} Suppose that $\ep X=0$ and $\ep X^2=\sigma^2<\infty$. Then
$$\lim_{\epsilon\searrow 0}\epsilon^2
\sum_{n=1}^{\infty}\frac 1{n\log n}
\pr(|S_n|\ge \epsilon\sqrt{ n\log\log n})=\sigma^2. $$
\endproclaim

The main purpose of this paper is
to show general results under the {\it minimal} conditions by using an
Feller's (1945) and Einmahl's (1989) truncation method.
The following two theorems are our main results.

\begin{theorem} \label{th1}
Let $a>-1$ and $b>-1/2$ and let $a_n(\epsilon)$ be a function of $\epsilon$ such that
\begin{eqnarray} \label{co1.1}
a_n(\epsilon) \log \log n \to \tau
\; \text{ as } \;
n\to \infty  \text{ and } \epsilon \searrow \sqrt{1+a}.
\end{eqnarray}
Suppose that
\begin{eqnarray} \label{co1.2}
\ep X=0, \; \ep X^2=\sigma^2<\infty
\; \text{ and } \;
\ep\big[ X^2(\log |X|)^a(\log\log |X|)^{b-1}\big]<\infty
\end{eqnarray}
and
\begin{eqnarray} \label{co1.3}
\ep X^2I\{|X|\ge t\}=o((\log\log t)^{-1}) \; \text{ as }\;
t\to \infty.
\end{eqnarray}
Then
\begin{eqnarray}\label{eqT1.1}
& & \lim_{\epsilon\searrow \sqrt{1+a}}
(\epsilon^2-a-1)^{b+1/2} \sum_{n=1}^{\infty}
\frac{(\log n)^a (\log\log n)^b}{n}
\pr \Big\{M_n \ge \sigma \phi(n) (\epsilon
+a_n(\epsilon))\Big\}
\non
& & \qquad \qquad \qquad=2\sqrt{\frac 1{\pi (a+1)} }\exp\{-2\tau\sqrt{1+a}\}
\Gamma(b+1/2)
\end{eqnarray}
and
\begin{eqnarray}\label{eqT1.2}
& & \lim_{\epsilon\searrow \sqrt{1+a}}
(\epsilon^2-a-1)^{b+1/2} \sum_{n=1}^{\infty}
\frac{(\log n)^a (\log\log n)^b}{n}
\pr \Big\{|S_n| \ge \sigma \phi(n) (\epsilon
+a_n(\epsilon))\Big\}
\non
& & \qquad \qquad \qquad=\sqrt{\frac 1{\pi (a+1)} }\exp\{-2\tau\sqrt{1+a}\}
\Gamma(b+1/2).
\end{eqnarray}
Here, $\Gamma(\cdot)$ is a gamma function.
Conversely, if (\ref{eqT1.1}) or (\ref{eqT1.2}) holds for $a>-1$, $b>-1/2$ and some $0<\sigma<\infty$, then (\ref{co1.2}) holds and
\begin{equation}\label{eqT1.3}
\liminf_{t\to \infty}(\log\log t)\ep X^2I\{|X|\ge t\}=0.
\end{equation}
\end{theorem}

\begin{theorem} \label{th2}
Suppose that $\ep X=0$ and $\ep X^2=\sigma^2<\infty$,
and let $a_n=O(1/\log\log n)$.
For $b>-1$, we have
\begin{eqnarray}\label{eqT2.1}
& & \lim_{\epsilon\searrow 0}
\epsilon^{2(b+1)} \sum_{n=1}^{\infty}
\frac{(\log\log n)^b}{n\log n} \pr \Big\{M_n \ge
\sigma\phi(n) (\epsilon+a_n)
\Big\}
\non
& & \qquad \qquad \qquad=
\frac 2{(b+1)\sqrt{\pi}}\Gamma(b+3/2)\sum_{k=0}^{\infty}
\frac{(-1)^k}{(2k+1)^{2b+2}}
\end{eqnarray}
and
\begin{eqnarray}\label{eqT2.2}
& & \lim_{\epsilon\searrow 0}
\epsilon^{2(b+1)} \sum_{n=1}^{\infty}
\frac{(\log\log n)^b}{n\log n} \pr \Big\{|S_n| \ge
\sigma\phi(n) (\epsilon+a_n)
\Big\}
\non
& & \qquad \qquad \qquad=
\frac 1{(b+1)\sqrt{\pi}}\Gamma(b+3/2).
\end{eqnarray}
Conversely, if (\ref{eqT2.1}) or (\ref{eqT2.2}) holds for some $b>-1$ and  $0<\sigma<\infty$, then $\ep X=0$ and $\ep X^2=\sigma^2$.
\end{theorem}

\begin{remark}
Note that the condition (\ref{co1.3}) is
sharp. A sufficient condition for it is given by
$$ \ep X^2 \log\log |X|<\infty. $$
So, (\ref{co1.3}) is weaker than Gut and Sp\u ataru's
condition in Theorem A (see also their Remark 1.1).
When $a>0$ (or $a=0$ and $b\ge 2$), the condition (\ref{co1.3}) is implied by
(\ref{co1.2}).
\end{remark}

\begin{remark} The condition that $\ep X=0$ and $\ep X^2<\infty$
is obviously sufficient and necessary for the conclusion of Theorem B
to hold,
by Theorem \ref{th2}. (see also Remark 1.2 of Gut and Sp\u ataru, 2000).
\end{remark}

The proofs of Theorem \ref{th1} and \ref{th2} are given in Section 4.
Before that,
we first verify (\ref{eqT1.1}), (\ref{eqT1.2}),
(\ref{eqT2.1}) and (\ref{eqT2.2}) under the assumption that $F$ is the normal distribution in Section 2,
after which, by using the truncation
and approximation method, we then show that the probabilities in
(\ref{eqT1.1}), (\ref{eqT1.2}),  (\ref{eqT2.1}) and
(\ref{eqT2.2}) can be replaced
by those for normal random variables in Section 3.
Throughout this
paper, we let $K(\alpha,\beta,\cdots)$, $C(\alpha,\beta,\cdots)$
etc denote positive constants which depend on $\alpha,\beta,
\cdots$ only, whose values can differ in different places.
$a_n\sim b_n$ means that $a_n/b_n\to 1$.

%%%%%%%%%%%%%%%%%%%%%%%%%%%%%%%%%%%%%%
%%%%%%%%%%%%%%%%%%%%%%%%%%%%%%%%%%%%%%%%

\section{Normal cases.}
\setcounter{equation}{0}

In this section, we prove Theorems \ref{th1} and \ref{th2} in the case that $\{X, X_n; n\ge 1\}$ are normal random variables.
Let $\{W(t); t\ge 0\}$ be a standard Wiener process and $N$ a standard normal variable. Our results are as follows.

\begin{proposition}\label{prop2.1}
Let $a>-1$ and $b>-1/2$ and let $a_n(\epsilon)$ be a function of $\epsilon$ satisfying (\ref{co1.1}).
Then
\begin{eqnarray}\label{eqprop2.1.1}
& & \lim_{\epsilon\searrow \sqrt{1+a}}
(\epsilon^2-a-1)^{b+1/2} \sum_{n=1}^{\infty}
\frac{(\log n)^a (\log\log n)^b}{n}
\non
& &  \qquad \qquad \qquad \qquad\cdot
\pr \Big\{\sup_{0\le s\le 1}|W(s)| \ge \sqrt{2\log\log n} (\epsilon
+a_n(\epsilon))\Big\}
\non
& & \qquad \qquad \qquad=2\sqrt{\frac 1{\pi (a+1)} }\exp\{-2\tau\sqrt{1+a}\}
\Gamma(b+1/2)
\end{eqnarray}
and
\begin{eqnarray}\label{eqprop2.1.2}
& & \lim_{\epsilon\searrow \sqrt{1+a}}
(\epsilon^2-a-1)^{b+1/2} \sum_{n=1}^{\infty}
\frac{(\log n)^a (\log\log n)^b}{n}
\pr \Big\{|N| \ge \sqrt{2\log\log n} (\epsilon
+a_n(\epsilon))\Big\}
\non
& & \qquad \qquad \qquad=\sqrt{\frac 1{\pi (a+1)} }\exp\{-2\tau\sqrt{1+a}\}
\Gamma(b+1/2)
\end{eqnarray}
\end{proposition}

\begin{proposition}\label{prop2.2} Let $a_n=O(1/\log\log n)$.
For any $b>-1$, we have
\begin{eqnarray*}
& & \lim_{\epsilon\searrow 0}
\epsilon^{2(b+1)} \sum_{n=1}^{\infty}
\frac{(\log\log n)^b}{n\log n} \pr \Big\{\sup_{0\le s\le 1}|W(s)|
 \ge
 (\epsilon+a_n)\sqrt{2\log\log n}
\Big\}
\\
& & \qquad \qquad \qquad=
\frac 2{(b+1)\sqrt{\pi}}\Gamma(b+3/2)\sum_{k=0}^{\infty}
\frac{(-1)^k}{(2k+1)^{2b+2}}
\end{eqnarray*}
and
\begin{eqnarray*}
& & \lim_{\epsilon\searrow 0}
\epsilon^{2(b+1)} \sum_{n=1}^{\infty}
\frac{(\log\log n)^b}{n\log n} \pr \Big\{|N| \ge
 (\epsilon +a_n)\sqrt{2\log\log n}
\Big\}
\\
& & \qquad \qquad \qquad=
\frac 1{(b+1)\sqrt{\pi}}\Gamma(b+3/2).
\end{eqnarray*}
\end{proposition}

The following lemma will be used in the proofs.

\begin{lemma}\label{lem2.1}
Let $\{W(t); t\ge 0\}$ be a standard Wiener process. Then for all $x>0$,
\begin{eqnarray}\label{eqL2.1.1}
 \pr \big(\sup_{0\le s\le 1}|W(s)| \ge x\big)
&=&
1-\sum_{k=-\infty}^{\infty}
(-1)^k\pr\big((2k-1)x\le N\le (2k+1)x\big)
\non
&=&4\sum_{k=0}^{\infty}
(-1)^k\pr\big( N\ge (2k+1)x\big)
\non
&=&2\sum_{k=0}^{\infty}
(-1)^k\pr\big( |N|\ge (2k+1)x\big).
\end{eqnarray}
In particular,
$$
\pr\big( \sup_{0\le s\le 1}|W(s)|\ge x \big)
\sim 2\pr\big( |N|\ge x\big)\sim \frac 4{\sqrt{2\pi}x}e^{-x^2/2}
\; \text{ as } \; x\to +\infty.
$$
\end{lemma}
{\bf Proof.} It is well known. See Billingsley (1968).

\bigskip

Now, we turn to prove the propositions.

{\noindent\bf Proof Proposition \ref{prop2.1}:} First, note that
the limit in (\ref{eqprop2.1.1}) does not depend on any finite
terms of the infinite series. Secondly, by Lemma \ref{lem2.1} and the condition (\ref{co1.1}) we
have
\begin{eqnarray*}
& &\pr \Big\{\sup_{0\le s\le 1}|W(s)|
   \ge \sqrt{2\log\log n}
(\epsilon +a_n(\epsilon))\Big\}
\sim 2\pr \Big\{|N|
   \ge \sqrt{2\log\log n}
(\epsilon +a_n(\epsilon))\Big\}
\\
& & \qquad \sim \frac 4{\sqrt{2\pi}(\epsilon+a_n(\epsilon))
 \sqrt{2\log\log n}}
\exp\Big\{-(\epsilon+a_n(\epsilon))^2 \log\log n \Big\} \\
& & \qquad \sim \frac 2{\sqrt{\pi}\epsilon }
\frac 1{\sqrt{\log\log n}}
\exp\Big\{-\epsilon^2\log\log n \Big\}
\exp\Big\{-2\epsilon a_n(\epsilon)\log\log n
\Big\}
\end{eqnarray*}
as $n\to \infty$, uniformly in
$\epsilon\in (\sqrt{1+a},\sqrt{1+a}+\delta)$
for some $\delta>0$. So, for any $0<\theta<1$, there exist $\delta>0$ and $n_0$ such that for all $n\ge n_0$ and
$\epsilon\in (\sqrt{1+a},\sqrt{1+a}+\delta)$,
$$\begin{array}{ll}
& \frac 2{\sqrt{\pi(1+a)} }\frac 1{\sqrt{\log\log n}}
\exp\Big\{-\epsilon^2\log\log n \Big\}
\exp\Big\{-2\tau\sqrt{1+a}-\theta\Big\}
\\
\le &\pr \Big\{\sup_{0\le s\le 1}|W(s)|
   \ge \sqrt{2\log\log n}
(\epsilon +a_n(\epsilon))\Big\}
\\
\le & \frac 2{\sqrt{\pi(1+a)} }\frac 1{\sqrt{\log\log n}}
\exp\Big\{-\epsilon^2\log\log n \Big\}
\exp\Big\{-2\tau\sqrt{1+a}+\theta\Big\}
\end{array} $$
and
$$\begin{array}{ll}
& \frac 1{\sqrt{\pi(1+a)} }\frac 1{\sqrt{\log\log n}}
\exp\Big\{-\epsilon^2\log\log n\Big\}
\exp\Big\{-2\tau\sqrt{1+a}-\theta\Big\}
\\
\le &\pr \Big\{|N|
   \ge \sqrt{2\log\log n}
(\epsilon +a_n(\epsilon))\Big\}
\\
\le & \frac 1{\sqrt{\pi(1+a)} }\frac 1{\sqrt{\log\log n}}
\exp\Big\{-\epsilon^2\log\log n \Big\}
\exp\Big\{-2\tau\sqrt{1+a}+\theta\Big\},
\end{array} $$
by the condition (\ref{co1.1}) again.
Also,
\begin{eqnarray*}
& &  \lim_{\epsilon\searrow \sqrt{1+a}}
(\epsilon^2-a-1)^{b+1/2}
\sum_{n=1}^{\infty}\frac{(\log n)^a (\log\log n)^b}{n}
\frac 1{\sqrt{\log\log n}}
\exp\Big\{-\epsilon^2 \log\log n\Big\} \\
&=& \lim_{\epsilon\searrow \sqrt{1+a} }
(\epsilon^2-a-1)^{b+1/2}
\int_{e^e}^{\infty} \frac{(\log x)^a (\log\log x)^{b-1/2}}{x}
\exp\Big\{-\epsilon^2\log\log x\Big\} dx \\
&=&
 \lim_{\epsilon\searrow \sqrt{1+a}}
(\epsilon^2-a-1)^{b+1/2}
\int_1^{\infty} y^{b-1/2}
\exp\Big\{-y\big(\epsilon^2-1-a\big)\Big\} dy \\
&=& \lim_{\epsilon\searrow \sqrt{1+a}}
\int_{\epsilon^2-1-a}^{\infty} y^{b-1/2}
e^{-y} dy
=  \int_0^{\infty}
y^{b-1/2}e^{-y} dy
 =\Gamma(b+1/2).
\end{eqnarray*}
(\ref{eqprop2.1.1}) and (\ref{eqprop2.1.2}) are now proved.

\bigskip

{\noindent\bf Proof Proposition \ref{prop2.2}:}
Observe  (\ref{eqL2.1.1}),
$$ \pr(|N|\ge x)=2\pr(N\ge x), \; \forall x>0, $$
and for any $m\ge 1$ and $x>0$,
\begin{eqnarray*}
 4\sum_{k=0}^{2m+1}
(-1)^k\pr\big( N\ge (2k+1)x\big)
&\le& \pr \big(\sup_{0\le s\le 1}|W(s)| \ge x\big)
\\
&\le&4\sum_{k=0}^{2m}
(-1)^k\pr\big( N\ge (2k+1)x\big).
\end{eqnarray*}
It is sufficient to show that for any $q>0$,
\begin{eqnarray*}
& & \lim_{\epsilon\searrow 0}
\epsilon^{2(b+1)} \sum_{n=1}^{\infty}
\frac{(\log\log n)^b}{n\log n}
\pr\Big(N\ge q(\epsilon+a_n)\sqrt{2\log\log n}\Big) \\
& & \qquad \qquad =q^{-2(b+1)}\frac 1{2(b+1)\sqrt{\pi}}\Gamma(b+3/2).
\end{eqnarray*}
Obviously,
\begin{eqnarray*}
& & \lim_{\epsilon\searrow 0}
\epsilon^{2(b+1)} \sum_{n=1}^{\infty}
\frac{(\log\log n)^b}{n\log n}
\pr\Big(N\ge q(\epsilon+a_n)\sqrt{2\log\log n}\Big)
\\
&=& q^{-2(b+1)}
 \lim_{\epsilon\searrow 0}
\epsilon^{2(b+1)} \sum_{n=1}^{\infty}
\frac{(\log\log n)^b}{n\log n}
\pr\Big(N\ge (\epsilon+q a_n)\sqrt{2\log\log n}\Big).
\end{eqnarray*}
So, it is sufficient to show that
\begin{eqnarray*}
 \lim_{\epsilon\searrow 0}
\epsilon^{2(b+1)} \sum_{n=1}^{\infty}
\frac{(\log\log n)^b}{n\log n}
\pr\Big(N\ge (\epsilon+a_n)\sqrt{2\log\log n}\Big)
=\frac 1{2(b+1)\sqrt{\pi}}\Gamma(b+3/2).
\end{eqnarray*}
Without losing of generality, we assume that $|a_n|\le \tau/\log\log n$.
Notice that
\begin{eqnarray*}
& & \Big|\pr\Big(N\ge (\epsilon+a_n)\sqrt{2\log\log n}\Big)
-\pr\Big(N\ge \epsilon\sqrt{2\log\log n}\Big)\Big|\\
&\le& \frac 1{\sqrt{2\pi}}
\exp\big\{-\frac{2\log\log n(\epsilon-\tau/\log\log n)^2}{2}\big\}
|a_n|\sqrt{2\log\log n} \\
&\le&\frac{\tau }{\sqrt{\log\log n}}
\exp\big\{-\epsilon^2\log\log n+2\epsilon\tau \big\}
 \end{eqnarray*}
and
\begin{eqnarray*}
 & &\lim_{\epsilon\searrow 0}
\epsilon^{2(b+1)} \sum_{n=1}^{\infty}
\frac{(\log\log n)^b}{n\log n}\cdot
\frac{1}{\sqrt{\log\log n}}
\exp\big\{-\epsilon^2\log\log n \big\} \\
&=& \lim_{\epsilon\searrow 0}
\epsilon^{2(b+1)} \int_{e^e}^{\infty}
\frac{(\log\log x)^{b-1/2}}{x\log x}
\exp\big\{-\epsilon^2\log\log x \big\} dx \\
&=& \lim_{\epsilon\searrow 0}
\epsilon^{2(b+1)} \int_1^{\infty}
y^{b-1/2}\exp\{-\epsilon^2 y \} dy
= \lim_{\epsilon\searrow 0}
\epsilon \int_{\epsilon^2}^{\infty}
y^{b-1/2}e^{-y} dy\\
&=&\lim_{\epsilon\searrow 0}
\epsilon \int_{\epsilon^2}^1
y^{b-1/2}e^{-y} dy
+\lim_{\epsilon\searrow 0}
\epsilon \int_1^{\infty}
y^{b-1/2}e^{-y} dy \\
&\le& \lim_{\epsilon\searrow 0}
\epsilon \int_{\epsilon^2}^1
y^{b-1/2} dy=0.
\end{eqnarray*}
Thus, it follows that
\begin{eqnarray*}
& & \lim_{\epsilon\searrow 0}
\epsilon^{2(b+1)} \sum_{n=1}^{\infty}
\frac{(\log\log n)^b}{n\log n}
\pr\Big(N\ge (\epsilon+a_n)\sqrt{2\log\log n}\Big)
\\
&=&\lim_{\epsilon\searrow 0}
\epsilon^{2(b+1)} \sum_{n=1}^{\infty}
\frac{(\log\log n)^b}{n\log n}
\pr\Big(N\ge \epsilon \sqrt{2\log\log n}\Big)\\
 &= &\lim_{\epsilon\searrow 0}
\epsilon^{2(b+1)} \int_{e^e}^{\infty}
\frac{(\log\log x)^b}{x\log x}
\pr\Big(N\ge \epsilon\sqrt{2\log\log x}\Big)dx
\\
&=&\lim_{\epsilon\searrow 0}
 \int_{\epsilon^2}^{\infty}
y^b\pr\Big(N\ge \sqrt{2y}\Big)dy
= \frac 1{b+1}\int_0^{\infty}
\pr\Big(N\ge \sqrt{2y}\Big)d y^{b+1} \\
&=&\frac 1{b+1}\pr\Big(N\ge \sqrt{2y}\Big) y^{b+1}\big|_0^{\infty}
+\frac 1{2(b+1)\sqrt{\pi} }\int_0^{\infty} y^{b+1/2}e^{-y}dy\\
&=&\frac 1{2(b+1)\sqrt{\pi}}\Gamma(b+3/2).
\end{eqnarray*}
The proposition is now proved.

%%%%%%%%%%%%%%%%%%%%%%%%%%%%%%%%%%%%%
%%%%%%%%%%%%%%%%%%%%%%%%%%%%%%%%%%%%

\section{Truncation and Approximation.}
\setcounter{equation}{0}

The purpose of this section is to use Feller's (1945)
 and Einmahl's (1989)
truncation methods to show that the probabilities in
(\ref{eqT1.1}), (\ref{eqT2.1}) for $M_n$ can be
approximated by those for $\sqrt{n}\sup_{0\le s\le 1}|W(s)|$
and the probabilities in
(\ref{eqT1.2}), (\ref{eqT2.2}) for $S_n$ can be
approximated by those for $\sqrt{n}N$.

Suppose that $\ep X=0$ and $\ep X^2=\sigma^2<\infty$.
Without losing of generality, we assume that $\sigma=1$ throughout this section.
 Let $p>1/2$. For each $n$ and $1\le j\le n$,
we let
\begin{eqnarray*}
X_{nj}^{\prime}=X_{nj}I\{|X_j|\le \sqrt{n}/(\log\log n)^p \},
& &\overline X_{nj}^{\prime}=X_{nj}^{\prime}-\ep [X_{nj}^{\prime}], \\
S_{nj}^{\prime}=\sum_{i=1}^jX_{nj}^{\prime},
& & \overline S_{nj}^{\prime}=\sum_{i=1}^j \overline X_{nj}^{\prime}, \\
\overline M_n^{\prime}=\max_{k\le n}|\overline S_{nk}^{\prime}|,
& &B_n=\sum_{k=1}^n \Var(\overline X_{nk}^{\prime})
\end{eqnarray*}
and
\begin{eqnarray*}
X_{nj}^{\prime\prime}=
X_{nj}I\{ \sqrt{n}/(\log\log n)^p<|X_j|\le \phi(n) \},
& &\overline X_{nj}^{\prime\prime}=X_{nj}^{\prime\prime}
  -\ep [X_{nj}^{\prime\prime}], \\
X_{nj}^{\prime\prime\prime}=
X_{nj}I\{ |X_j|> \phi(n) \},
& &\overline X_{nj}^{\prime\prime\prime}
=X_{nj}^{\prime\prime\prime}-\ep [X_{nj}^{\prime\prime\prime}].
\end{eqnarray*}
And also define $S_{nj}^{\prime\prime}$,
$S_{nj}^{\prime\prime\prime}$,
$\overline S_{nj}^{\prime\prime}$,
$\overline S_{nj}^{\prime\prime\prime}$,
$\overline M_n^{\prime\prime}$ and
$\overline M_n^{\prime\prime\prime}$ similarly.

The following two propositions are the main results of this section.

\begin{proposition}\label{prop3.1}
Let $a>-1$, $b>-1$ and $2\ge p>p^{\prime}>1/2$.
Suppose that the condition (\ref{co1.2}) is satisfied.
Then there exist $\delta>0$ and a sequence of positive numbers $\{q_n \}$ such that
\begin{eqnarray} \label{eqprop3.1.1}
&&\pr\Big\{\sup_{0\le s\le 1}|W(s)|
   \ge \epsilon \sqrt{2\log\log n }
 + \frac{3}{ (\log\log n)^{p^{\prime}} } \Big\}
-q_n
\non &\le&
  \pr\Big\{M_n \ge \epsilon \sqrt{2B_n\log\log n } \Big\}
\non
&\le&\pr\Big\{\sup_{0\le s\le 1}|W(s)|
     \ge \epsilon \sqrt{2\log\log n }
  - \frac{3}{ (\log\log n)^{p^{\prime}} } \Big\} +q_n,
\end{eqnarray}
\begin{eqnarray} \label{eqprop3.1.1b}
&&\pr\Big\{|N| \ge \epsilon \sqrt{2\log\log n }
 + \frac{3}{ (\log\log n)^{p^{\prime}} } \Big\}
-q_n
\non &\le&
 \pr\Big\{|S_n| \ge \epsilon \sqrt{2B_n\log\log n } \Big\}
\non
&\le&\pr\Big\{|N|   \ge \epsilon
    \sqrt{2\log\log n }
   -\frac{3}{(\log\log n)^{p^{\prime}} } \Big\} +q_n,
\non & &  \qquad \forall
\epsilon\in (\sqrt{1+a}-\delta, \sqrt{1+a}+\delta),
\quad n\ge 1
\end{eqnarray}
and
\begin{eqnarray}\label{eqprop3.1.2}
\sum_{n=1}^{\infty}\frac{(\log n)^a(\log\log n)^b}{n}
q_n \le K(a,b,p,p^{\prime},\delta)<\infty.
\end{eqnarray}
\end{proposition}

\begin{proposition}\label{prop3.2}
Let $b$ be a real number and $2\ge p>p^{\prime}>1/2$.
Suppose that $\ep X=0$ and $\ep X^2=1$. Then
\begin{eqnarray}\label{eqprop3.2.1}
& &
\pr\big(\sup_{0\le s\le 1}|W(s)|\ge x+3/(\log\log n)^{p^{\prime}}\big)-q_n^{\ast}
\le \pr\big(M_n\ge x\sqrt{B_n}\big)
\non
&\le&
\pr\big(\sup_{0\le s\le 1}|W(s)|\ge x-3/(\log\log n)^{p^{\prime}}\big)+q_n^{\ast},
\quad
\forall x>0,
\end{eqnarray}
\begin{eqnarray}\label{eqprop3.2.1b}
& &
\pr\big(|N|\ge x+3/(\log\log n)^{p^{\prime}}\big)-q_n^{\ast}
\le \pr\big(|S_n|\ge x\sqrt{B_n}\big)
\non
&\le&
\pr\big(|N|\ge x-3/(\log\log n)^{p^{\prime}}\big)+q_n^{\ast},
\quad
\forall x>0,
\end{eqnarray}
where $q_n^{\ast}\ge 0$ satisfies
\begin{equation} \label{eqprop3.2.2}
\sum_{n=1}^{\infty} \frac{(\log\log n)^b}{n\log n} q_n^{\ast} \le
K(b,p, p^{\prime})<\infty.
\end{equation}
\end{proposition}

To show this two results, we need some lemmas.

\begin{lemma}\label{lem1}
For any sequence of independent random variables $\{\xi_n; n\ge
1\}$ with mean zero and finite variance, there exists a sequence
of independent normal variables $\{\eta_n; n\ge 1\}$ with $\ep
\eta_n=0$ and $\ep \eta_n^2=\ep \xi_n^2$ such that, for all $Q>2$
and $y>0$,
$$ \pr\Big(\max_{k\le n}|\sum_{i=1}^k \xi_i-\sum_{i=1}^k \eta_i|\ge y\Big)
\le (AQ)^Qy^{-Q}\sum_{i=1}^n \ep |\xi_i|^Q, $$
 whenever $\ep|\xi_i|^Q<\infty$, $i=1,\ldots,n$.
Here, $A$ is a universal
constant.
\end{lemma}

{\bf Proof.} See Sakhaneko (1980,1984, 1985).

\begin{lemma}\label{lem2} Let $Q\ge 2$, $\xi_1, \xi_2,\ldots, \xi_n$ be independent random variables with $\ep \xi_k=0$ and $\ep |\xi_k|^Q<\infty$,
$k=1,\ldots, n$. Then for all $y>0$,
$$ \pr\Big(\max_{k\le n}|\sum_{i=1}^k \xi_i|\ge y \Big)
\le 2\exp\Big\{-\frac{y^2}{8\sum_{k=1}^n \Var(\xi_k) }\Big\}
+(2AQ)^Qy^{-Q}\sum_{i=1}^n \ep |\xi_i|^Q, $$
where $A$ is a universal constant as in Lemma \ref{lem1}.
\end{lemma}

{\bf Proof.} It follows from Lemma \ref{lem1} easily. See also Petrov (1995, Page 78).

\begin{lemma}\label{lem3}
Define $\Delta_n=\max_{k\le n}|\overline S_{nk}^{\prime}-S_k|$.
Let $a>-1$, $b>-1$ and $p>1/2$. Suppose that the condition (\ref{co1.2}) is satisfied and $\ep X^2=1$.
Then for any $\lambda>0$ there exist a constant $K=K(a,b,p,\lambda)$ such that
\begin{eqnarray}\label{eqL3.1}
\sum_{n=1}^{\infty} \frac{(\log n)^a(\log\log n)^b}{n}I_n
\le K \ep \Big[X^2(\log |X|)^a(\log\log |X|)^{b-1} \Big]
 <\infty,
\end{eqnarray}
where
$$
I_n=\pr\Big(\Delta_n\ge \sqrt{n}/(\log\log n)^2,
 \overline M_n^{\prime}\ge \lambda \phi(n)\Big).
$$
\end{lemma}
{\bf Proof.}
Let $\beta_n=n \ep[|X|I\{|X|>\sqrt{n}/(\log\log n)^p\}]$. Then
$|\ep\sum_{i=1}^jX_{ni}^{\prime}|\le \beta_n$, $1\le j\le n$.
Setting
$$\Cal L=\{n:\beta_n\le \frac 18 \sqrt{n}/(\log\log n)^2\},$$
we have
$$\{\Delta_n\ge \sqrt{n}/(\log\log n)^2\}
\subset
\bigcup_{j=1}^n \{ X_j\ne X_{nj}^{\prime}\}, \quad n\in \Cal L. $$
So for $n\in \Cal L$,
\begin{eqnarray*}
I_n&\le& \sum_{j=1}^n \pr\Big(X_j\ne X_{nj}^{\prime}, \overline M_n^{\prime}\ge \lambda\phi(n) \Big).
\end{eqnarray*}
Observer that $X_{nj}^{\prime}=0$ whenever $X_j\ne X_{nj}^{\prime}$,
$j\le n$, so that  we have for $n$ large enough and all $1\le j\le n$,
\begin{eqnarray*}
& &\pr\Big(X_j\ne X_{nj}^{\prime},
\overline M_n^{\prime}\ge \lambda\phi(n) \Big)\\
&=&\pr\Big(X_j\ne X_{nj}^{\prime},
\max_{k\le j-1}|\overline S_{nk}^{\prime}|\vee
\max_{j<k\le n}|\overline S_{nk}^{\prime}-X_{nj}^{\prime}|
\ge \lambda\phi(n) \Big)\\
&=&\pr\Big(X_j\ne X_{nj}^{\prime}\Big)
\pr\Big(\max_{k\le j-1}|\overline S_{nk}^{\prime}|\vee
\max_{j<k\le n}|\overline S_{nk}^{\prime}-X_{nj}^{\prime}|
\ge \lambda\phi(n) \Big)\\
&\le &\pr\Big(X_j\ne X_{nj}^{\prime}\Big)
\pr\Big(\overline M_n^{\prime}\ge \lambda\phi(n)
-|X_{nj}^{\prime}|\Big)\\
&\le &\pr\Big(|X|>\sqrt{n}/(\log\log n)^p \Big)
\pr\Big(\overline M_n^{\prime}\ge \lambda\phi(n)
-\sqrt{n}/(\log\log n)^p\Big)\\
&\le &\pr\Big(|X|>\sqrt{n}/(\log\log n)^p \Big)
\pr\Big(\overline M_n^{\prime}\ge \frac{\lambda}{2}\phi(n)
\Big).
\end{eqnarray*}
A straightforward application of the inequalities of Ottaviani and Bernstein yields:
\begin{eqnarray*}
 \pr\Big(\overline M_n^{\prime}\ge \frac{\lambda}{2}\phi(n)
\Big)
&\le& 2\pr\Big(|\overline S_n^{\prime}|\ge \frac{\lambda}{4}\phi(n)
\Big)
\le (\log n)^{-\eta} \\
& & \quad \text{ for some } \eta=\eta(\lambda)>0.
\end{eqnarray*}
 So,
\begin{eqnarray*}
& & \sum_{n\in \Cal L} \frac{(\log n)^a (\log\log n)^b }{n}
I_n \\
&\le& C\sum_{n=1}^{\infty}\frac{(\log n)^a (\log\log n)^b }{n}
\cdot n \pr\Big(|X|>\frac{\sqrt{n}}{(\log\log n)^p} \Big)(\log n)^{-\eta}\\
&\le&  \sum_{n=1}^{\infty}\sum_{j=n}^{\infty}
 \pr\Big(\frac{\sqrt{j}}{(\log \log j)^p}<|X|
  \le \frac{\sqrt{j+1}}{(\log\log (j+1))^p}\Big)
(\log n)^{a-\eta}(\log\log n)^b \\
&\le&  \sum_{j=1}^{\infty}
 \pr\Big(\frac{\sqrt{j}}{(\log \log j)^p}<|X|
  \le \frac{\sqrt{j+1}}{(\log\log (j+1))^p}\Big)
\sum_{n=1}^j (\log n)^{a-\eta }(\log\log n)^b \\
&\le&  \sum_{j=1}^{\infty}
 \pr\Big(\frac{\sqrt{j}}{(\log \log j)^p}<|X|
  \le \frac{\sqrt{j+1}}{(\log\log (j+1))^p}\Big)
j (\log j)^{a-\eta }(\log\log j)^b \\
&\le& C\ep\Big[X^2(\log |X|)^{a-\eta }(\log\log |X|)^{b+2p} \Big]
\le C\ep \Big[X^2(\log |X|)^a(\log\log |X|)^{b-1} \Big].
\end{eqnarray*}
If $n\not\in \Cal L$, then we have
\begin{eqnarray*}
 I_n
\le \pr\Big(\overline M_n^{\prime}
  \ge  \lambda \phi(n) \Big)
\le  (\log n)^{-\eta}.
\end{eqnarray*}
It follows that
\begin{eqnarray*}
& & \sum_{n\not\in \Cal L}
\frac{(\log n)^a(\log\log n)^b}{n}
I_n
\le \sum_{n\not\in \Cal L}
\frac{(\log n)^{a-\eta}(\log\log n)^b}{n} \\
&\le& 8\sum_{n\not\in \Cal L}
\frac{(\log n)^{a-\eta}(\log\log n)^{b+2}}{n^{3/2}}
 \beta_n  \\
&\le& 8\sum_{n=1}^{\infty} \frac{(\log n)^{a-\eta}(\log\log n)^{b+2}}{n^{1/2}}
\\
& & \qquad \cdot \sum_{j=n}^{\infty}
\ep\Big[|X|I\big\{\frac{\sqrt{j}}{(\log \log j)^p}<|X|
  \le \frac{\sqrt{j+1}}{(\log\log (j+1))^p}\big\}\Big]\\
&=& 8\sum_{j=1}^{\infty}
\ep\Big[|X|I\big\{\frac{\sqrt{j}}{(\log \log j)^p}<|X|
  \le \frac{\sqrt{j+1}}{(\log\log (j+1))^p}\big\}\Big] \\
& & \qquad \cdot \sum_{n=1}^j
\frac{(\log n)^{a-\eta}(\log\log n)^{b+2}}{n^{1/2}}\\
&\le& C\sum_{j=1}^{\infty}
 \ep\Big[|X|I\big\{\frac{\sqrt{j}}{(\log \log j)^p}<|X|
  \le \frac{\sqrt{j+1}}{(\log\log (j+1))^p}\big\}\Big]\\
& & \qquad \cdot
\sqrt{j}(\log j)^{a-\eta}(\log\log j)^{b+2}\\
&\le& C\ep\Big[X^2(\log |X|)^{a-\eta}
(\log\log |X|)^{b+2+p}\Big]
\le C\ep \Big[X^2(\log |X|)^a(\log\log |X|)^{b-1} \Big].
\end{eqnarray*}
(\ref{eqL3.1}) is proved.

\begin{lemma}\label{lem4}
Let $a>-1$, $b>-1$ and $p>1/2$. Suppose the condition (\ref{co1.2}) is satisfied and $\ep X^2=1$.
Then for any $\lambda>0$ there exist a constant $K=K(a,b,p,\lambda)$ such that
$$\sum_{n=1}^{\infty} \frac{(\log n)^a(\log\log n)^b}{n}II_n
\le K \ep \Big[X^2(\log |X|)^a(\log\log |X|)^{b-1} \Big]
 <\infty, $$
where
$$II_n=\pr\Big(\Delta_n\ge \sqrt{n}/(\log\log n)^2,
 M_n \ge \lambda\phi(n)\Big). $$
\end{lemma}
{\bf Proof.} Obviously,
$$ II_n\le \pr\Big(\Delta_n\ge \sqrt{n}/(\log\log n)^2,
 \overline M_n^{\prime} \ge \frac{\lambda}{3}\phi(n)\Big)
+\pr\Big(\overline M_n^{\prime\prime} \ge \frac{\lambda}{3}\phi(n)\Big)
+\pr\Big(\overline M_n^{\prime\prime\prime} \ge \frac{\lambda}{3}\phi(n)\Big).
$$
Observe that
$\max_{k\le n}|\ep S_{nk}^{\prime\prime\prime}|
\le n \ep X^2 /\phi(n)=o(\sqrt n)$. We have
\begin{eqnarray*}
& & \sum_{n=1}^{\infty} \frac{(\log n)^a(\log\log n)^b}{n}
\pr\Big(\overline M_n^{\prime\prime\prime} \ge \frac{\lambda}{3}\phi(n)\Big) \\
&\le& C\sum_{n=1}^{\infty} \frac{(\log n)^a(\log\log n)^b}{n}
\sum_{j=1}^n \pr\Big(X_j^{\prime\prime\prime}\ne 0\Big)\\
&\le&\sum_{n=1}^{\infty} (\log n)^a (\log\log n)^b\pr\big(|X|\ge \phi(n)\big) \\
&\le& K \ep \Big[X^2(\log |X|)^a(\log\log |X|)^{b-1} \Big].
\end{eqnarray*}
Also, notice that
$\sum_{k=1}^n \Var(\overline X_{nk}^{\prime\prime})\le
n \ep \big[X^2I\big\{\sqrt{n}/(\log\log n)^p<|X|\le \phi(n)\big\}\big]
=o(n)$. By Lemma \ref{lem2} we have for $Q>2$,
\begin{eqnarray*}
& & \sum_{n=1}^{\infty} \frac{(\log n)^a(\log\log n)^b}{n}
\pr\Big(\overline M_n^{\prime\prime} \ge \frac{\lambda}{3}\phi(n)\Big) \\
&\le& C \sum_{n=1}^{\infty} \frac{(\log n)^a(\log\log n)^b}{n}
\exp\Big\{ - \frac{\lambda^2\phi^2(n)}{3^2 8\cdot o(n)}\Big\} \\
& & \quad + C \sum_{n=1}^{\infty} \frac{(\log n)^a(\log\log n)^b}{n}
\cdot \frac {3^Q}{\lambda^Q \phi^Q(n)}n \ep\big[|X|^QI\{|X|\le \phi(n)\}\big] \\
&\le& K+C\sum_{n=1}^{\infty} \frac{(\log n)^a(\log\log n)^b}{\phi^Q(n)}
\sum_{j=1}^n\ep\big[|X|^QI\{\phi(j-1)<|X|\le \phi(j)\}\big] \\
&\le& K+C\sum_{j=1}^{\infty}\ep\big[|X|^QI\{\phi(j-1)<|X|\le \phi(j)\}\big] \sum_{n=j}^{\infty} \frac{(\log n)^a(\log\log n)^b}{\phi^Q(n)}
\\
&\le& K+C\sum_{j=1}^{\infty}\ep\big[|X|^QI\{\phi(j-1)<|X|\le \phi(j)\}\big]  j\frac{(\log j)^a(\log\log j)^b}{\phi^Q(j)}\\
&\le& K+C\sum_{j=1}^{\infty}\ep\big[|X|^2I\{\phi(j-1)<|X|\le \phi(j)\}\big]  (\log j)^a(\log\log j)^{b-1} \\
&\le& K+C \ep \Big[X^2(\log |X|)^a(\log\log |X|)^{b-1} \Big]<\infty.
\end{eqnarray*}
Finally, by noticing Lemma \ref{lem3}, we compete the proof of Lemma \ref{lem4}.

\begin{lemma}\label{lem6} Suppose that the condition (\ref{co1.2}) is satisfied. Then for any $1/2<p^{\prime}<p$ we have
\begin{eqnarray}\label{eqL6.1}
& &\pr\big(\sup_{0\le s\le 1}|W(s)|
  \ge x+1/(\log\log n)^{p^{\prime}}\big)-p_n
\le\pr\big(\overline M_n^{\prime}\ge x\sqrt{B_n}\big)
\non
&\le&
 \pr\big(\sup_{0\le s\le 1}|W(s)|
 \ge x- 1/(\log\log n)^{p^{\prime}}\big)+p_n,
\quad \forall x>0
\end{eqnarray}
and
\begin{eqnarray}\label{eqL6.1b}
& &\pr\big(|N|\ge x+1/(\log\log n)^{p^{\prime}}\big)-p_n
\le\pr\big(|\overline S_n^{\prime}|\ge x\sqrt{B_n}\big)
\non
&\le&
 \pr\big(|N|\ge x-1/(\log\log n)^{p^{\prime}} \big)+p_n,
\quad \forall x>0,
\end{eqnarray}
 where $p_n\ge 0$ satisfies
\begin{equation}\label{eqL6.2}
\sum_{n=1}^{\infty}\frac{(\log n)^a(\log\log n)^b}{n} p_n
\le K(a,b,p,p^{\prime})<\infty.
\end{equation}
\end{lemma}

{\bf Proof.} By Lemma \ref{lem1}, there exist a universal constant
$A>0$ and a sequence of standard Wiener processes $\{W_n(\cdot)\}$
such that for all $Q>2$,
\begin{eqnarray*}
& & \pr\Big( \max_{k\le n}|\overline S_{nk}^{\prime}-W_n(\frac kn B_n)|
\ge \frac 12 \sqrt{B_n}/(\log\log n)^{p^{\prime}} \Big) \\
&\le& (AQ)^Q\Big(\frac{ (\log\log n)^{p^{\prime}} }{\sqrt{B_n}}\Big)^Q
\sum_{k=1}^n \ep\big|\overline X_{nk}^{\prime}\big|^Q \\
&\le&  C n \Big(\frac{(\log\log n)^{p^{\prime}}}{\sqrt{n}}\Big)^Q
 \ep\big[|X|^QI\{|X|\le \sqrt{n}/(\log \log n)^p\}\big].
\end{eqnarray*}
On the other hand, by Lemma 1.1.1 of Cs\"org\H o and R\'ev\'esz (1981),
\begin{eqnarray*}
& & \pr\Big(|\max_{0\le s\le 1}|W_n(s B_n)-W_n(\frac{[ns]}{n} B_n)|
\ge \frac 12 \sqrt{B_n}/(\log\log n)^{ p^{\prime} }\Big) \\
&=& \pr\Big(\max_{0\le s\le 1}|W_n(s)-W_n(\frac{[ns]}{n} )|
\ge \frac 12 \sqrt{\frac 1n}
\frac{\sqrt n}{ (\log\log n)^{p^{\prime}} }\Big) \\
&\le& Cn
\exp\Big\{
-\frac{(\frac 12 \sqrt{n}/(\log\log n)^{p^{\prime}} )^2}{3}\Big\}
\le C n\exp\Big\{-\frac 1{12}n/(\log\log n)^{2p^{\prime}}\Big\}.
\end{eqnarray*}
Let
\begin{equation} \label{eqL6.3}
p_n=
\pr\Big(\sup_{0\le s\le 1}
\Big| \frac{\overline S_{n,[ns]}^{\prime}}{\sqrt{B_n}}
-\frac{W_n(sB_n)}{\sqrt{B_n}}\Big|
\ge \frac 1{(\log\log n)^{p^{\prime}} } \Big).
\end{equation}
Then $p_n$ satisfies (\ref{eqL6.1}) and (\ref{eqL6.1b}), since $\{W_n(t
B_n)/\sqrt{B_n}; t\ge 0\}\overset{\Cal D}=\{W(t); t\ge 0\}$ for
each $n$. And also,
$$p_n\le C n \Big(\frac{ (\log\log n)^{p^{\prime}} }{\sqrt{n}}\Big)^Q
 \ep\big[|X|^QI\{|X|\le \sqrt{n}/(\log\log n)^p \}\big]
+C n\exp\Big\{-\frac 1{12}n/(\log\log n)^{2p^{\prime}}\Big\}. $$
It follows that
\begin{eqnarray*}
& &\sum_{n=1}^{\infty}\frac{(\log n)^a(\log\log n)^b}{n} p_n \\
&\le& K_1+
C\sum_{n=1}^{\infty}
\frac{(\log n)^a(\log\log n)^{ b+ p^{\prime}Q }  }{n^{Q/2}}
\ep\big[|X|^QI\{|X|\le \sqrt{n}/(\log\log n)^p\}\big] \\
&\le& K_1+
C\sum_{n=1}^{\infty}
\frac{(\log n)^a(\log\log n)^{b+ p^{\prime}Q } }{n^{Q/2}}\\
& & \quad \cdot
\sum_{j=1}^n\ep\big[|X|^Q
 I\big\{\frac{\sqrt{j-1}}{(\log\log(j-1) )^p}<|X|\le
     \frac{\sqrt{j}}{(\log\log j)^p}\big\} \big] \\
&\le& K_1+
C\sum_{j=1}^{\infty}
\ep\big[|X|^Q
  I\big\{\frac{\sqrt{j-1}}{(\log\log(j-1) )^p}<|X|\le
     \frac{\sqrt{j}}{(\log\log j)^p}\big\} \big]
j\frac{(\log j)^a(\log\log j)^{  b+ p^{\prime}Q }  }{j^{Q/2}}
\\
&\le& K_1+
C\ep\big[|X|^2(\log|X|)^a(\log\log |X|)^{b+ (p^{\prime}-p)Q +2p}\big]
\le K<\infty,
\end{eqnarray*}
whenever $(p^{\prime}-p)Q +2p<-1$. So, (\ref{eqL6.2}) is satisfied.

%%%%%%%%%%%%%%%%%%%%%%%%%%
\bigskip
Now, we turn to prove Propositions \ref{prop3.1} and \ref{prop3.2}.

{\noindent\bf Proof of Proposition \ref{prop3.1}:}
 Let $0<\delta<\frac 14\sqrt{1+a}$.
Observe that, if $n$ is large enough,
\begin{eqnarray*}
& & \pr\Big\{M_n \ge \epsilon \sqrt{2 B_n\log\log n}
 \Big\}\\
&=&\pr\Big\{M_n \ge \epsilon
\sqrt{2 B_n\log\log n},
\Delta_n\le \frac{\sqrt{n}}{(\log\log n)^2} \Big\} \\
& & \quad  +\pr\Big\{M_n \ge \epsilon \sqrt{2 B_n\log\log n},
\Delta_n> \frac{\sqrt{n}}{(\log\log n)^2} \Big\}\\
&\le&\pr\Big\{\overline M_n^{\prime} \ge \epsilon \sqrt{2 B_n\log\log n}
- \frac{\sqrt{n}}{(\log\log n)^2} \Big\} \\
& & \quad  +\pr\Big\{M_n \ge \frac {\sqrt{1+a}}{4} \phi(n),
\Delta_n> \frac{\sqrt{n}}{(\log\log n)^2} \Big\}\\
&\le&\pr\Big\{\overline M_n^{\prime}
\ge \sqrt{B_n}\big[\epsilon \sqrt{2 \log\log n}- \frac 2{(\log\log n)^2}\big] \Big\} +II_n \\
&\le&\pr\Big\{\sup_{0\le s\le 1}|W(s)| \ge
\epsilon \sqrt{2 \log\log n}- \frac {2}{(\log\log n)^2}
-\frac{1}{(\log\log n)^{p^{\prime}} }
\Big\} +p_n +II_n \\
&\le&\pr\Big\{\sup_{0\le s\le 1}|W(s)| \ge
\epsilon \sqrt{2 \log\log n}
-\frac{3}{(\log\log n)^{p^{\prime}} }
\Big\} +p_n +II_n
\end{eqnarray*}
for all $\epsilon\in(\sqrt{1+a}-\delta,\sqrt{1+a}+\delta)$,
where $II_n$ is defined in Lemmas \ref{lem4}
with $\lambda=\sqrt{1+a}/4$ and $p_n$ is defined in \ref{lem6}.
Also, if  $n$ is large enough,
\begin{eqnarray*}
& & \pr\Big\{M_n \ge \epsilon \sqrt{2 B_n\log\log n}
\Big\}\\
&\ge&\pr\Big\{M_n \ge \epsilon \sqrt{2 B_n\log\log n},
\Delta_n\le \frac{\sqrt{n}}{(\log\log n)^2} \Big\} \\
&\ge&\pr\Big\{\overline M_n^{\prime}
  \ge \epsilon \sqrt{2 B_n\log\log n}
 + \frac{\sqrt{n}}{(\log\log n)^2},
\Delta_n\le \frac{\sqrt{n}}{(\log\log n)^2}  \Big\} \\
&\ge&\pr\Big\{\overline M_n^{\prime} \ge \sqrt{B_n} \big[\epsilon
\sqrt{2\log\log n}  + \frac{2}{(\log\log n)^2}\big]
\Big\}
\\
& &\quad -\pr\Big\{\overline M_n^{\prime}
 \ge \frac {\sqrt{1+a}}{4} \phi(n),
 \Delta_n> \frac{\sqrt{n}}{(\log\log n)^2}\Big\} \\
&\ge&\pr\Big\{\sup_{0\le s\le 1}|W(s)| \ge \epsilon
\sqrt{2\log\log n}  + \frac{3}{(\log\log n)^{p^{\prime}} }\Big\}
-p_n -I_n
\end{eqnarray*}
for all $\epsilon\in(\sqrt{1+a}-\delta,\sqrt{1+a}+\delta)$,
where $I_n$ is  defined in Lemma \ref{lem3}
 with $\lambda=\sqrt{1+a}/4$.

Similarly,  if  $n$ is large enough,
\begin{eqnarray*}
&&\pr\Big\{|N| \ge \epsilon \sqrt{2\log\log n }
 + \frac{3}{ (\log\log n)^{p^{\prime}} } \Big\}
-p_n-I_n
\\
&\le&
 \pr\Big\{|S_n| \ge \epsilon \sqrt{2B_n\log\log n } \Big\}
\\
&\le&\pr\Big\{|N|   \ge \epsilon
    \sqrt{2\log\log n }
   -\frac{3}{(\log\log n)^{p^{\prime}} } \Big\} +p_n+II_n.
\end{eqnarray*}
Letting $q_n=p_n+I_n+II_n$ completes the proof by
Lemmas \ref{lem3}, \ref{lem4} and \ref{lem6}.

\bigskip

%%%%%%%%%%%%%%%%%%%%%%%%%
{\noindent\bf Proof Proposition \ref{prop3.2}:} Let
$\{W_n(\cdot)\}$ be a sequence of standard Wiener processes being
defined in the proof of Lemma \ref{lem6}, and let $p_n$ be defined
in (\ref{eqL6.3}). And set
$$q_n^{\ast}=
\pr\Big(\sup_{0\le s\le 1}\big|M_{[ns]}/\sqrt{B_n}-W_n(sB_n)/\sqrt{B_n}\big|\ge 3/(\log\log n)^{p^{\prime}}\Big). $$
Then $q_n^{\ast}$ satisfies (\ref{eqprop3.2.1}) and (\ref{eqprop3.2.1b}),
 and also
$$q_n^{\ast}\le \pr\big(\Delta_n\ge \sqrt{n}/(\log\log n)^2\big)+p_n. $$
By Lemma \ref{lem6},
$$\sum_{n=1}^{\infty}
\frac{(\log\log n)^b}{n\log n} p_n \le K_1(b,p, p^{\prime})<\infty.
$$
Also, following the lines in the proof of (\ref{eqL3.1}) we have
\begin{eqnarray*}
& &\sum_{n=1}^{\infty}\frac{(\log\log n)^b}{n\log n}
   \pr\big(\Delta_n\ge \sqrt{n}/(\log\log n)^2\big) \\
&\le&\sum_{n\in \Cal L}\frac{(\log\log n)^b}{n\log n}\cdot n
   \pr\big(|X|> \sqrt{n}/(\log \log n)^p \big)
+\sum_{n\not\in \Cal L}\frac{(\log\log n)^{b+2}}{n^{3/2}\log n}
  \beta_n\\
&\le&\sum_{n=1}^{\infty}\frac{(\log\log n)^b}{\log n}
   \pr\big(|X|> \sqrt{n}/(\log \log n)^p \big) \\
& & \quad
+\sum_{n=1}^{\infty}\frac{(\log\log n)^{b+2}}{\sqrt{n}\log n}
   \ep\big[|X|I\{|X|> \sqrt{n}/(\log \log n)^p \}\big]
\\
&\le&C\ep\big[X^2(\log|X|)^{-1}(\log\log |X|)^{b+2p}\big]
  +C\ep\big[X^2(\log|X|)^{-1}(\log\log |X|)^{b+2+p}\big]\\
&\le& C\ep X^2<\infty.
\end{eqnarray*}
So, $q_n^{\ast}$ satisfies (\ref{eqprop3.2.2}).

%%%%%%%%%%%%%%%%%%%%%%%%%%%%%%%%%%%%%%%%%%%%%%
%%%%%%%%%%%%%%%%%%%%%%%%%%%%%%%%%%%%%%%%%%%%

\section{Proofs of the Theorems.}
\setcounter{equation}{0}

%%%%%%%%%%%%%%%
\subsection{Proofs of the direct parts.}

Without losing of generality, we assume that $\ep X=0$ and $\ep X^2=1$.

{\noindent\bf Proof of the direct part of Theorem \ref{th1}:}
Let $\delta>0$ small enough and $\{q_n\}$ be such that
(\ref{eqprop3.1.1}), (\ref{eqprop3.1.1b}) and (\ref{eqprop3.1.2}) hold.
Then
\begin{eqnarray*}
 \lim_{\epsilon\searrow \sqrt{1+a}}
(\epsilon^2-a-1)^{b+1/2} \sum_{n=1}^{\infty}
\frac{(\log n)^a (\log\log n)^b}{n} q_n=0,
\end{eqnarray*}
by (\ref{eqprop3.1.2}). Notice that $a_n(\epsilon)\to 0$.
By (\ref{eqprop3.1.1}), we have that for $n$ large enough,
\begin{eqnarray*}
&&\pr\Big\{\sup_{0\le s\le 1}|W(s)| \ge
\sqrt{2\log\log n}
(\epsilon+a_n(\epsilon))
 + \frac{3}{(\log\log n)^{p^{\prime}} } \Big\}
-q_n
\\
&\le& \pr\Big\{M_n \ge
\sqrt{2 B_n \log\log n}
(\epsilon+a_n(\epsilon)) \Big\}
\\
&\le&\pr\Big\{\sup_{0\le s\le 1}|W(s)| \ge
\sqrt{2\log\log n}
(\epsilon+a_n(\epsilon))
- \frac{3}{ (\log\log n)^{p^{\prime}} } \Big\}
+q_n,
\\
& &  \qquad \forall \epsilon\in
(\sqrt{1+a}-\delta/2, \sqrt{1+a}+\delta/2).
\end{eqnarray*}
On the other hand, by Proposition \ref{prop2.1},
\begin{eqnarray*}
& & \lim_{\epsilon\searrow \sqrt{1+a}}
(\epsilon^2-a-1)^{b+1/2} \sum_{n=1}^{\infty}
\frac{(\log n)^a (\log\log n)^b}{n}
\non
& &  \qquad \qquad \qquad \qquad\cdot
\pr \Big\{\sup_{0\le s\le 1}|W(s)|
 \ge \sqrt{2\log\log n} (\epsilon
+a_n(\epsilon))
\pm \frac{3}{(\log\log n)^{p^{\prime}} }\Big\}
\\
& & \qquad \qquad \qquad=2\sqrt{\frac 1{\pi (a+1)} }\exp\{-2\tau\sqrt{1+a}\}
\Gamma(b+1/2).
\end{eqnarray*}
It follows that
\begin{eqnarray}\label{eqP1.1}
& & \lim_{\epsilon\searrow \sqrt{1+a}}
(\epsilon^2-a-1)^{b+1/2} \sum_{n=1}^{\infty}
\frac{(\log n)^a (\log\log n)^b}{n}
\non
& &  \qquad \qquad \qquad \qquad\cdot
\pr \Big\{M_n \ge \sqrt{2B_n \log\log n} (\epsilon
+a_n(\epsilon))\Big\}
\non
& & \qquad \qquad \qquad=2\sqrt{\frac 1{\pi (a+1)} }\exp\{-2\tau\sqrt{1+a}\}
\Gamma(b+1/2).
\end{eqnarray}
Similarly,
\begin{eqnarray}\label{eqP1.1b}
& & \lim_{\epsilon\searrow \sqrt{1+a}}
(\epsilon^2-a-1)^{b+1/2} \sum_{n=1}^{\infty}
\frac{(\log n)^a (\log\log n)^b}{n}
\non
& &  \qquad \qquad \qquad \qquad\cdot
\pr \Big\{|S_n| \ge \sqrt{2B_n \log\log n} (\epsilon
+a_n(\epsilon))\Big\}
\non
& & \qquad \qquad \qquad=\sqrt{\frac 1{\pi (a+1)} }\exp\{-2\tau\sqrt{1+a}\}
\Gamma(b+1/2).
\end{eqnarray}

Finally, noticing the condition (\ref{co1.3}), we have
$$0\le n-B_n\le 2 n \ep[X^2I\{|X|\ge \sqrt{n}/(\log\log n)^p \}]
=o(n(\log\log n)^{-1}). $$
Let $a_n^{\prime}(\epsilon)=\sqrt{n/B_n}(\epsilon+a_n(\epsilon))-\epsilon$.
Then
$$\pr\Big\{M_n \ge \phi(n)(\epsilon +a_n(\epsilon))\Big\}
=\pr\Big\{M_n \ge
\sqrt{2 B_n\log\log n}
(\epsilon +a_n^{\prime}(\epsilon))\Big\}, $$
$$\pr\Big\{|S_n| \ge \phi(n)(\epsilon +a_n(\epsilon))\Big\}
=\pr\Big\{ |S_n| \ge
\sqrt{2 B_n\log\log n}
(\epsilon +a_n^{\prime}(\epsilon))\Big\}, $$
and,
$$a_n^{\prime}(\epsilon)\log\log n
=\epsilon \frac{(n-B_n)\log\log n}{\sqrt{B_n}(\sqrt n+\sqrt{B_n})}+\sqrt{\frac{n}{B_n}}a_n(\epsilon)\log\log n
\to \tau $$
as $n\to \infty$ and $\epsilon\searrow \sqrt{1+a}$.
Now, (\ref{eqT1.1}) and (\ref{eqT1.2}) follow from (\ref{eqP1.1})
and (\ref{eqP1.1b}), respectively.

%%%%%%%%%%%%%%%

\bigskip

{\noindent \bf Proof of the direct part of Theorem \ref{th2}:}
We show (\ref{eqT2.1}) only, since the proof of (\ref{eqT2.2}) is similar.
Noticing $n\ge B_n\sim n$ and Proposition \ref{prop3.2}, for any $0<\delta<1$ we have for $n$ large enough and all $\epsilon>0$,
\begin{eqnarray*}
& & \pr\Big\{\sup_{0\le s\le 1}|W(s)|
\ge [\epsilon(1+\delta)+2|a_n|+3/\log\log n]
\sqrt{2\log\log n}\Big\}-q_n^{\ast} \\
&\le& \pr\Big\{\sup_{0\le s\le 1}|W(s)|\ge (\epsilon+a_n)(1+\delta)
\sqrt{2\log\log n} + 3/(\log\log n)^{p^{\prime}}\Big\}
  -q_n^{\ast} \\
&\le&\pr\Big\{M_n\ge (\epsilon +a_n)(1+\delta)\sqrt{2 B_n \log\log n} \Big\}\\
&\le&\pr\Big\{M_n\ge (\epsilon+a_n)\phi(n)\Big\}
\le \pr\Big\{M_n\ge
 (\epsilon+a_n)
 \sqrt{2 B_n \log\log n}
\Big\} \\
&\le&\pr\Big\{\sup_{0\le s\le 1}|W(s)|\ge
(\epsilon+a_n) \sqrt{2\log\log n}
  - 3/(\log\log n)^{p^{\prime}} \Big\}
  +q_n^{\ast} \\
&\le&\pr\Big\{\sup_{0\le s\le 1}|W(s)|\ge (\epsilon+a_n-3/\log\log n)
\sqrt{2\log\log n} \Big\}
  +q_n^{\ast}.
\end{eqnarray*}
So, by Propositions \ref{prop2.2} and \ref{prop3.2},
\begin{eqnarray*}
& &(1+\delta)^{-2(b+1)}
\frac 2{(b+1)\sqrt{\pi}}\Gamma(b+3/2)\sum_{k=0}^{\infty}
\frac{(-1)^k}{(2k+1)^{2b+2}}
\\
&\le& \liminf_{\epsilon\searrow 0}
\epsilon^{2(b+1)} \sum_{n=1}^{\infty}
\frac{(\log\log n)^b}{n\log n}
\pr\Big\{M_n\ge (\epsilon+a_n)\phi(n) \Big\} \\
&\le& \limsup_{\epsilon\searrow 0}
\epsilon^{2(b+1)} \sum_{n=1}^{\infty}
\frac{(\log\log n)^b}{n\log n}
\pr\Big\{M_n\ge (\epsilon+a_n)\phi(n) \Big\} \\
&\le& \frac 2{(b+1)\sqrt{\pi}}\Gamma(b+3/2)\sum_{k=0}^{\infty}
\frac{(-1)^k}{(2k+1)^{2b+2}}.
\end{eqnarray*}
Letting $\delta\to 0$ completes the proof.

%%%%%%%%%%%%%%%
\bigskip

%%%%%%%%%%%%%%%%%%%%
%%%%%%%%%%%

\subsection{Proofs of the converse parts.}
Now, we turn to prove the converse parts of Theorem
\ref{th1} and \ref{th2}. First, we show that each of (\ref{eqT1.1}),
(\ref{eqT1.2}), (\ref{eqT2.1}) and (\ref{eqT2.2}) implies
\begin{eqnarray} \label{eqpf4.3}
\ep X^2<\infty,\; \; \ep X=0 \;
\; \text{ and } \;
\ep\big[ X^2(\log |X|)^a(\log\log |X|)^{b-1}\big]<\infty,
\end{eqnarray}
where $a=-1$ in Theorem \ref{th2}.
We only give the proof that (\ref{eqT1.2}) implies (\ref{eqpf4.3}),
since other proofs are similar.
Let $\{\widetilde X, \widetilde X_n; n\ge 1\}$ be the symmetrization of $\{X, X_n; n\ge 1\}$, and let
$\widetilde S_n=\sum_{k=1}^n\widetilde X_k$.
Then by (\ref{eqT1.2}),
\begin{eqnarray*}
 \limsup_{\epsilon\searrow \sqrt{1+a}}
(\epsilon^2-a-1)^{b+1/2} \sum_{n=1}^{\infty}
\frac{(\log n)^a (\log\log n)^b}{n}
\pr \Big\{|\widetilde S_n| \ge 2\sigma \phi(n) (\epsilon
+a_n(\epsilon))\Big\}
\le K<\infty.
\end{eqnarray*}
For $M>0$, define $Y=Y(M)=\widetilde XI\{|\widetilde X|<M\}$ and
$Y_n=Y_n(M)=\widetilde X_nI\{|\widetilde X_n|<M\}$.
Observing that
$\widetilde XI\{|\widetilde X|<M\}-\widetilde XI\{|\widetilde X|\ge M\}\overset{\Cal D}=\widetilde X$ and
$\widetilde XI\{|\widetilde X|<M\}-\widetilde XI\{|\widetilde X|\ge M\}+\widetilde X=2 Y$, we obtain that
\begin{eqnarray}\label{eqpf4.4}
& & \limsup_{\epsilon\searrow \sqrt{1+a}}
(\epsilon^2-a-1)^{b+1/2} \sum_{n=1}^{\infty}
\frac{(\log n)^a (\log\log n)^b}{n}
\pr \Big\{|\sum_{k=1}^n Y_k| \ge 2\sigma \phi(n) (\epsilon
+a_n(\epsilon))\Big\}
\non
& &\le
 2\limsup_{\epsilon\searrow \sqrt{1+a}}
(\epsilon^2-a-1)^{b+1/2} \sum_{n=1}^{\infty}
\frac{(\log n)^a (\log\log n)^b}{n}
\pr \Big\{|\widetilde S_n| \ge 2\sigma \phi(n) (\epsilon
+a_n(\epsilon))\Big\}
\non
& &\le 2K<\infty.
\end{eqnarray}
However, since $Y$ is a bounded random variable which satisfies
conditions (\ref{co1.2}) and (\ref{co1.3}),
by the direct part of Theorem \ref{th1} we have
\begin{eqnarray}\label{eqpf4.5}
& & \lim_{\epsilon\searrow \sqrt{1+a}}
(\epsilon^2-a-1)^{b+1/2} \sum_{n=1}^{\infty}
\frac{(\log n)^a (\log\log n)^b}{n}
\pr \Big\{|\sum_{k=1}^n Y_k| \ge \sqrt{\ep Y^2} \phi(n) (\epsilon
+a_n(\epsilon))\Big\}
\non
& & \qquad \qquad \qquad=\sqrt{\frac 1{\pi (a+1)} }\exp\{-2\tau\sqrt{1+a}\}
\Gamma(b+1/2)>0.
\end{eqnarray}
Putting (\ref{eqpf4.4}) and (\ref{eqpf4.5}) together yields
$\sqrt{\ep \widetilde X^2I\{|\widetilde X|<M\}}=\sqrt{\ep Y^2}\le 2\sigma$.
Then, letting $M\to \infty$ yields $\ep X^2<\infty$.

$\ep X=0$ is obvious when $\ep X^2<\infty$, for otherwise we have
$$\pr\{|S_n|\ge \epsilon \sigma \phi(n)\}\to 1, \quad \forall \epsilon>0,  $$
which implies that
$$\sum_{n=1}^{\infty}
\frac{(\log n)^a (\log\log n)^b}{n}
\pr \{|S_n| \ge \epsilon\sigma \phi(n) \}=\infty,
\quad \forall \epsilon>0,\; a\ge -1 \; \text{ and } \; b\ge -1. $$

Now,
by (\ref{eqT1.2}) and the L\'evy inequality we obtain that for some $\epsilon>0$,
\begin{eqnarray*}
& &\sum_{n=1}^{\infty}
\frac{(\log n)^a (\log\log n)^b}{n}
\pr\{\max_{k\le n}|X_k|\ge 3\epsilon \sigma \phi(n) \} \\
& &\le C\sum_{n=1}^{\infty}
\frac{(\log n)^a (\log\log n)^b}{n}
\pr\{\max_{k\le n}|X_k|\ge 2\epsilon \sigma \phi(n)
+2\sqrt{n\ep X^2} \} \\
& & \le C\sum_{n=1}^{\infty}
\frac{(\log n)^a (\log\log n)^b}{n}
\pr \{\max_{k\le n}|S_k| \ge\epsilon \sigma \phi(n)+\sqrt{n\ep X^2} \} \\
& & \le C\sum_{n=1}^{\infty}
\frac{(\log n)^a (\log\log n)^b}{n}
\pr \{|S_n| \ge\epsilon \sigma \phi(n) \} <\infty.
\end{eqnarray*}
Observe that
$$\pr\{\max_{k\le n}|X_k|\ge 3\epsilon \sigma \phi(n) \}
\le \frac{\ep X^2}{18\epsilon^2\log\log n}\to 0. $$
We conclude that
\begin{eqnarray*}
& &\sum_{n=1}^{\infty}
(\log n)^a (\log\log n)^b
\pr\{|X|\ge 3\epsilon \sigma \phi(n) \} \\
&\le &C\sum_{n=1}^{\infty}
\frac{(\log n)^a (\log\log n)^b}{n}
\pr\{\max_{k\le n}|X_k|\ge 3\epsilon \sigma \phi(n) \}
<\infty,
\end{eqnarray*}
which implies
$$ \ep\big[ X^2(\log |X|)^a(\log\log |X|)^{b-1}\big]<\infty. $$
(\ref{eqpf4.3}) is proved.

 Next, we show that $\ep X^2=\sigma^2$.
 By the direct part of Theorem
\ref{th2}, (\ref{eqT2.1}) and (\ref{eqT2.2}) shall hold with $\ep X^2$ taking the
place of $\sigma^2$, which are obviously contradictory to (\ref{eqT2.1}) and (\ref{eqT2.2}) themselves, respectively,  if $\ep
X^2\ne \sigma^2$. Notice that (\ref{eqP1.1}) and (\ref{eqP1.1b})
hold whenever (\ref{eqpf4.3}) is satisfied. However, if $\ep
X^2\ne \sigma^2$, (\ref{eqT1.1}) and (\ref{eqT1.2}) are
contradictory to (\ref{eqP1.1}) and (\ref{eqP1.1b}), respectively,
since $B_n\sim n \ep X^2$.

Finally, we show (\ref{eqT1.3}). Suppose that
(\ref{eqT1.3}) fails. Without losing of generality, we can assume
that $\sigma^{-2}\ep [X^2I\{|X|\ge \sqrt{n}/(\log\log n)^p\}]\ge
\tau_0/\log\log n$ for some $\tau_0>0$ and all $n\ge 1$. Then
$n\sigma^2-B_n\ge n\ep [X^2I\{|X|\ge \sqrt{n}/(\log\log n)^p\}]\ge
n\sigma^2\tau_0/\log\log n$. Let
$a_n^{\prime}(\epsilon)=\sqrt{1+\tau_0/\log\log
n}\big(\epsilon+a_n(\epsilon)\big)-\epsilon$. Then
$$ a_n^{\prime}(\epsilon)\log\log n\to \tau+\tau_0\sqrt{1+a}/2,$$
and
$$\pr\Big\{M_n\ge \sigma\phi(n)\big(\epsilon+a_n(\epsilon)\big)\Big\}
 \le \pr\Big\{M_n\ge \sqrt{2B_n\log\log n}
 \big(\epsilon+a_n^{\prime}(\epsilon)\big)\Big\}, $$
$$\pr\Big\{|S_n|\ge \sigma\phi(n)\big(\epsilon+a_n(\epsilon)\big)\Big\}
 \le \pr\Big\{|S_n|\ge \sqrt{2B_n\log\log n}
 \big(\epsilon+a_n^{\prime}(\epsilon)\big)\Big\}, $$
It follows that (\ref{eqT1.1}) and (\ref{eqT1.2}) are
contradictory to (\ref{eqP1.1}) and (\ref{eqP1.1b}), respectively.
The proof is now completed.
%%%%%%%%%%%%%%%%
\newpage

%%%%%%%%%%%%%%%%%%%%%%%%%%%%%%%%%%%%%%%%%%%%%
%\begin{center}{\bf References}\end{center}


\begin{thebibliography}{99}
\bibitem{} Billingsley, P. (1968),
{\em Convergence of Probability Measures}.
Wiley, New York.

\bibitem{} Cs\"org\H o, M. and R\'ev\'esz, P.  (1981),
{\em Strong Approximations in Probability and Statistics}.
 Academic, New York.

\bibitem{} Einmahl, U. (1989),
The Darling-Erd\"os theorem for sums of i.i.d. random variables.
{\em Probab. Theory Relat. Fields}, {\bf 82}: 241-257.


\bibitem{} Feller, W. (1945),
The law of the iterated logarithm for identically distributed random variables.
{\em Ann. Math.}, {\bf 47}: 631-638.

\bibitem{} Gut, A. and Sp\u ataru, A. (2000),
Precise asymptotics in the law of the iterated logarithm.
{\em Ann. Probab.}, {\bf 28}:1870-1883.

\bibitem{} Jain, N. C. and Pruitt, W. E. (1975),
The other law of the iterated logarithm.
{\em Ann. Probab.}, {\bf 3}:1046-1049.

\bibitem{} Petrov, V. V. (1995),
{\em Limit Theorems of Probability Theory}.
Oxford University Press, Oxford.

\bibitem{} Sakhanenko, A. I.  (1980),
On unimprovable estimates of the rate of convergence in the invariance principle.
 In {\em Colloquia Math. Soci. J\'anos Bolyai}, {\bf 32}, pp. 779-783,
Nonparametric Statistical Inference, Budapest (Hungary).

\bibitem{} Sakhanenko, A. I. (1984),
On estimates of the rate of convergence in the invariance principle.
 In {\em Advances in Probab. Theory: Limit Theorems and Related Problems}
(A.A. Borovkov, Ed.),  pp. 124-135, Springer, New York.

\bibitem{} Sakhanenko, A. I. (1985),
Convergence rate in  the invariance principle for non-identically distributed variables with exponential moments.
In {\em Advances in Probab. Theory: Limit Theorems for Sums of Random Variables}
(A.A. Borovkov, Ed.),  pp. 2-73, Springer, New York.


\end{thebibliography}
\end{document}